# THE EXISTENCE OF A MOORE GRAPH OF DEGREE 57 IS STILL OPEN

*Vance Faber*
*Jonathan Keegan*
*revision: February 27, 2023*

*Note: This version is a complete rewrite incorporating ideas of the second author to find a more complete set of solutions.*

**Abstract**. In 2020, a paper [1] appeared in the arXiv claiming to prove that a Moore graph of diameter 2 and degree 57 does not exist. (The paper is in Russian – we include a link [4] to a translation of this paper kindly provided to us by Konstantin Selivanov.) The proof technique is reasonable. It employs the fact that from [2] such a graph must be distance regular and from [3] that there exists a large set of relations which such a graph must satisfy. The argument proceeds by a case analysis that shows that this set of relations cannot be satisfied. We show that this seems not to be correct. The system of equations has a coefficient matrix which decomposes into independent blocks of reasonable size, all of which have non-negative integer solutions. In addition, we show there are other inequality constraints besides non-negativity that a Moore graph imposes on solutions and enumerate all the solutions which satisfy these constraints. As an alternative to this method, we show that there is a family of systems of permutations with the property that the Moore graph exists if and only if there is a member of the family with no solutions.

**Intersection arrays**. In this section, we define some symbols and graphs related to the proof. We assume there is a Moore graph $G$ with degree $57$. Then by Theorem 7 in [2], $G$ has a distance-regular subgraph $\Gamma$ with intersection array $[55, 54, 2; 1, 1, 54]$. This graph is constructed as follows. Let $x$ and $y$ be any two adjacent vertices in $G$. Let $G_2(x)$ be the set of vertices $u$ where $d(u, x) = 2$. Then $\Gamma$ is the subgraph of $G$ induced by the vertices in $G_2(x) \cap G_2(y)$. What it means for a diameter $3$ distance-regular graph $\Gamma$ to have intersection array $[55, 54, 2; 1, 1, 54]$ is best illustrated on triples $(x, y, z)$. We fix $x$ and $y$ and count the number of neighbors $z$ of $y$ that satisfy $d(x, y) = Z$, $d(x, z) = Y$. Then these six numbers in order are given by

1) if $Y = 1, Z = 0$ the number is $55$ (the degree of the graph);
2) if $Y = 2, Z = 1$ is $54$;
3) if $Y = 3, Z = 2$ the number is $2$;
4) if $Y = 0, Z = 1$ the number is $1$;
5) if $Y = 1, Z = 2$ the number is $1$;
6) if $Y = 2, Z = 3$ the number is $54$.

We calculate three auxiliary symmetric matrices from this. Again, we appeal to the triples $(x, y, z)$ but this time, we allow $X = d(y, z)$ to vary. For fixed $Z$, $p_{XY}^Z$ is a $3 \times 3$ symmetric matrix which counts the number of vertices $z$ of distance $X$ from $y$ and distance $Y$ from $x$ where $Z = d(x, y)$ with $X, Y, Z \in \{1, 2, 3\}$. From the intersection array, we are given $p_{21}^1 = 54$, $p_{31}^2 = 2$, $p_{11}^2 = 1$, $p_{21}^3 = 54$. Apparently, this is enough information to compute all three of the intersection number matrices $p_{XY}^Z$ [1]. (This is also proved in Section 2.6 of https://arxiv.org/abs/1410.6294.)

The only thing we need now comes from the beginning of the paper [1]. It takes any triple of vertices, $u, v, w$ and counts the number of vertices $z$ which have distance $zu = i, zv = j, zw = k$. We assign the variable $x(i, j, k, u, v, w)$ to this number. Since the diameter of the distance-regular subgraph is 3, these distances are all between 0 and 3 (0 if $z$ is one of the vertices.). We also let $W = d(u, v)$, $U = d(v, w)$, $V = d(u, w)$ which are values between 1 and 3. Here I have to correct a misunderstanding in an earlier version pointed out to me by the second author. Given $i, j, k$, $x(i, j, k, u, v, w)$ is not dependent only on $U, V, W$. The equations that the $x$ satisfy depend only on $U, V, W$ but they have multiple solutions and one triple of vertices with distances $U, V, W$ might belong to one solution and a different triple of vertices with distances $U, V, W$ might belong to a different solution. However, the equations for $x$ depend only on $U, V, W$ so we may as well index the variables as $x(i, j, k, U, V, W)$. Note that there is symmetry which the paper mentions. So given $u, v, w$ and distances from $z$ to them, some variables can be deduced from others in the cases where some distances between the given vertices are equal (see symmetry below – the squiggles in [1] just denote various permutations of the three entries. The possible values for the variables are solutions to the linear system (+) below constrained by the fact that the entries must be non-negative integers and additional constraints on the variables discussed below. When a 0 is present as one of $i, j, k$ then this next paragraph gives the fixed value.

Let $u, v, w$ be vertices of a graph $\Gamma$ and $W = d(u, v)$, $U = d(v, w)$, $V = d(u, w)$. There is only one vertex $z = u$ such that $d(z, u) = 0$. Both the paper [3] and [1] say that therefore $x(0, j, h, u, v, w)$ is equal to 0 or 1 relative to the fixed triple $u, v, w$ depending on whether or not the distances $W = d(u, v)$ and $V = d(u, w)$ agree with $j$ and $h$, respectively. Let $\delta_{ij}$ be the Kronecker function. Then this says that
$x(0, j, h, u, v, w) = \delta_{jW} \delta_{hV}$, $x(i, 0, h, u, v, w) = \delta_{iW} \delta_{hU}$ and $x(i, j, 0, u, v, w) = \delta_{iU} \delta_{jV}$. I believe this is bit misleading. I believe that, for example, $x(0, j, h, u, v, w)$ depends on $U$ also. An example is $x(0, 1, 1, u, v, w)$ with $U = 3$. In this case $x(0, 1, 1, u, v, w) = 0$ even if $W = j = 1$ and $V = h = 1$. The problem occurs because a triple with distances $1,1,3$ fails the triangle inequality.

Now to derive the equations (+) discussed in [1], we fix the distance between any two vertices from $\{u, v, w\}$ and count the number of vertices of any distance to the third one:

$$\sum_{l=1}^{d} x(l, j, h, U, V, W) = p_{jh}^{U} - x(0, j, h, U, V, W)$$

$$\sum_{l=1}^{d} x(i, l, h, U, V, W) = p_{ih}^{V} - x(i, 0, h, U, V, W)$$

$$\sum_{l=1}^{d} x(i, j, l, U, V, W) = p_{ij}^{W} - x(i, j, 0, U, V, W).$$

Some of the triplets disappear. When $|i - j| > W$ or $|i + j| < W$, we have $p_{ij}^{W} = 0$ and then $x(i, j, h, U, V, W) = 0$ for all $0 \leq h \leq 3$.

**Equations and right-hand sides.** From the paper [1; Lemma 1] (plus the girth, the triangle inequality and triangle symmetry)

$$p^1 = \begin{pmatrix} 0 & 54 & 0 \\ 54 & 2808 & 108 \\ 0 & 108 & 2 \end{pmatrix},$$

$$p^2 = \begin{pmatrix} 1 & 52 & 2 \\ 54 & 2811 & 106 \\ 2 & 106 & 2 \end{pmatrix},$$

$$p^3 = \begin{pmatrix} 0 & 54 & 1 \\ 54 & 2862 & 54 \\ 1 & 54 & 54 \end{pmatrix},$$

We find it easier to relabel the indices $i, j, h, U, V, W$ by $i_1, i_2, i_3, i_4, i_5, i_6$. In this notation, the equations in (+) become

$$\sum_{i_1=0}^{d} x(i_1 i_2 i_3 i_4 i_5 i_6) = p_{i_2 i_3}^{i_4}, \quad \sum_{i_2=0}^{d} x(i_1 i_2 i_3 i_4 i_5 i_6) = p_{i_1 i_3}^{i_5}, \quad \sum_{i_3=0}^{d} x(i_1 i_2 i_3 i_4 i_5 i_6) = p_{i_1 i_2}^{i_6}.$$

One thing to note is that for fixed $U, V, W$, the equations that have variables with $(i_4, i_5, i_6) = (U, V, W)$ are completely independent from all the other equations. We collect these equations into a single system that we call a block. We label the block $\{U, V, W\}$ with $UVW$. Although the equations for different blocks are independent, it may be the constraints are not. We will discuss this further below. Because the entries with index 0 are known, we can make those values part of the right hand side. This gives us a system of equations

$$\sum_{i_1=1}^{d} x(i_1 i_2 i_3 i_4 i_5 i_6) = b(i_2 i_3 i_4 i_5 i_6) = \varepsilon(i_4 i_5 i_6)(p_{i_2 i_3}^{i_4} - \delta_{i_2 i_6} \delta_{i_3 i_5})$$

$$\sum_{i_2=1}^{d} x(i_1 i_2 i_3 i_4 i_5 i_6) = b(i_1 i_3 i_4 i_5 i_6) = \varepsilon(i_4 i_5 i_6)(p_{i_1 i_3}^{i_5} - \delta_{i_4 i_6} \delta_{i_1 i_3})$$

$$\sum_{i_1=1}^{d} x(i_1 i_2 i_3 i_4 i_5 i_6) = b(i_2 i_3 i_4 i_5 i_6) = \varepsilon(i_4 i_5 i_6)(p_{i_1 i_2}^{i_6} - \delta_{i_2 i_4} \delta_{i_1 i_5})$$

where $1 \leq i_j \leq 3$ and $\varepsilon(i_4 i_5 i_6) = 1$ if the triple $i_4 i_5 i_6$ is allowed and $0$ otherwise. Also, any of the other $i_j$ could take on the value $0$ but my calculations show that leads to tautologies so I have omitted them.

An important thing to note is that the coefficient matrices for all the blocks are identical. We call this matrix $M$. The problem is a linearly constrained integer programming problem.

**Null space of the coefficient matrix.** To describe all the solutions to a block, we start with a particular solution and add members of the null space while satisfying various constraints. To write the equations in standard form, we order the 27 variables in each block lexicographically: $i_1 i_2 i_3 \leq j_1 j_2 j_3$ if $i_1 \leq j_1$ and if $i_1 = j_1$ then $i_2 \leq j_2$ and if $i_1 i_2 = j_1 j_2$ then $i_3 \leq j_3$. For example, 133 is variable 9 while 331 is variable 25.

We fix three vertices $p, q, r$. The linear system for (+) can be written as $M\alpha = \beta$ where the matrix $M$ is completely independent of $p, q, r$. If we use the lexicographical order for the vertices and equations within a block, one way to write M is

$$M = \begin{bmatrix} [111] \otimes I_3 \otimes I_3 \\ I_3 \otimes [111] \otimes I_3 \\ I_3 \otimes I_3 \otimes [111] \end{bmatrix}.$$

To find the null space, we write the null vector as

$$X = \begin{bmatrix} u_1 \\ u_2 \\ u_3 \\ v_1 \\ v_2 \\ v_3 \\ w_1 \\ w_2 \\ w_3 \end{bmatrix}$$

where each of these entries is a vector in 3 dimensions. We also use the notation $\bar{u} = [111]u$, that is three times the mean of the entries of $u$. Then $MX = 0$ tells us that all the $u_i, v_i, w_i$ have mean zero and the following matrix has all row and column sums equal to zero:

$$\begin{bmatrix} u_1 & u_2 & u_3 \\ v_1 & v_2 & v_3 \\ w_1 & w_2 & w_3 \end{bmatrix}.$$

Any set of mean zero vectors $u_1, u_2, v_1, v_2$ uniquely determines the remaining 5 mean zero vectors using

$$u_3 = -u_1 - u_2$$
$$v_3 = -v_1 - v_2$$
$$w_1 = -u_1 - v_1$$
$$w_2 = -u_2 - v_2$$
$$w_3 = -w_1 - w_2 = -u_3 - v_3$$

This shows that the null space of $M$ has dimension 8 and a set of basis vectors is given by choosing one of $u_1, u_2, v_1, v_2$ to be either

$$\varepsilon_1 = \begin{bmatrix} 1 \\ 0 \\ -1 \end{bmatrix} \text{ or } \varepsilon_2 = \begin{bmatrix} 0 \\ 1 \\ -1 \end{bmatrix}$$

and the rest to be 0. These solutions can be written as $X_{ijk} = \varepsilon_i \otimes \varepsilon_j \otimes \varepsilon_k$ with $i, j \in \{1, 2\}$. For example

$$X_{111} = \begin{bmatrix} \varepsilon_1 \\ 0 \\ -\varepsilon_1 \\ 0 \\ 0 \\ 0 \\ -\varepsilon_1 \\ 0 \\ \varepsilon_1 \end{bmatrix}.$$

**Symmetries**. To understand the symmetries of the system, we think of the tuple as a triple of pairs $((i_1,i_4),(i_2,i_5),(i_3,i_6))$. Every permutation of the three pairs just relabels the configuration $u,v,w$ while fixing $z$. Because of the symmetry of the indices and the right-hand sides, the permutation $(i_1i_2)(i_4i_5)$ maps the first family of equations to the second and the permutation $(i_1i_3)(i_4i_6)$ maps the first family to the third. This means that the system of equations is invariant under the action of the 6 permutations $\sigma$ in the symmetric group $S_3$. Thus given a solution vector $x(i,j,k,U,V,W)$ to the block $UVW$ associated with the vertices $u,v,w$, the vector $x(\sigma(i,j,k),\sigma(vw,uw,uv))$ is a solution to the block determined by the vertices $\sigma(u,v,w)$. This does not necessarily mean however, that if $uw=uv$ then $x(i,j,k,vw,uw,uv) = x(i,k,j,vw,uv,uw)$. Some solutions we give below illustrate how that seemingly reasonable property fails. For example, there is a solution in block 322 where $x(24) = x(323322) \neq x(332322) = x(26)$.

**Known constraints and canonical blocks.** Some of these variables are known to have the value zero. Some entire blocks have only the trivial solution. By the triangle inequality, no three sides of a triangle can be $113$. This means that none of these triples can be $113, 131$ or $311$: $i_1i_2i_6$, $i_1i_3i_5$, $i_2i_3i_4$ and $i_4i_5i_6$. In addition, our graph has girth 5 so no triangle of distances can add up to $3$ and no square with all sides of length 1 can exist. This means that no triple $i_1i_2i_6$, $i_1i_3i_5$, $i_2i_3i_4$ or $i_4i_5i_6$ can be 111 and none of $i_1i_2i_4i_5$, $i_1i_3i_4i_6$, $i_2i_3i_5i_6$ can be 1111. Note that $112, 121$ and $211$ are allowed because the triple might be a path between the three vertices.

These considerations allow us to choose a set of non-trivial canonical blocks whose solutions determine all solutions by symmetry. We choose as our canonical set the 8 blocks 211, 221, 222, 321, 322, 331, 332 and 333.

**Additional constraints**. There are two important sets of additional constraints. In order prove them, we have to consider some properties of the graph $\Gamma = <G_2(x) \cap G_2(y)>$ where $x$ and $y$ are two adjacent vertices in the Moore graph $G$. Let $\Gamma_3$ be the graph consisting of the vertices of $\Gamma$ with an edge between vertices if the distance between them is 3.

*Theorem 1 (see, for example, [2]).* $\Gamma_3$ is isomorphic to the Cartesian product graph $K_{56} \otimes K_{56}$.

*Proof.* Because $G$ is a Moore graph, there are 56 vertices in $G$ which are adjacent to $x$ but not in the neighborhood of $y$. Each of these vertices $x_i$ has a set $X_i$ of 56 neighbors which are distance 2 from $x$ and distance 3 from $y$. So each of these 56 sets $X_i$ is a complete graph in $\Gamma_3$. But the same is true if we start from $y$ and get the 56 complete graphs $Y_j$. We just need to show that the intersection of $X_i$ and $Y_j$ is unique for each pair $i$ and $j$. Consider a $y_j$ which is a neighbor of $y$ and distance 2 from $x$. It cannot be adjacent to two members of an $X_i$ or it will be in a square with $x_i$. So each $Y_j$ intersects an $X_i$ at most once. But by the pigeon hole principle, every $Y_j$ must intersect each $X_i$ exactly once. This proves the theorem.

In what follows, we call the $X_i$ and $Y_j$ *grid lines*. Two vertices are distance 3 apart if and only if they are on a grid line. Each vertex is on exactly two grid lines, an $X_i$ and $Y_j$. Given any two non-intersecting grid lines, there is a matching between them in $\Gamma$.

*Lemma 2.* Given any $U \geq V \geq W$, the value of $x(3,3,3,U,V,W)$ is fixed:

a) if none of $U,V,W$ equals 3 then $x(3,3,3,U,V,W) = 0$;
b) if exactly one of $U,V,W$ equals 3 then $x(3,3,3,U,V,W) = 1$;
c) if exactly two of $U,V,W$ equals 3 then $x(3,3,3,U,V,W) = 0$.
d) if $U = V = W = 3$ then $x(3,3,3,U,V,W) = 53$.

*Proof.* We utilize the grid lines. We have

a) if none of $U,V,W$ equals 3 then $\Gamma_3(u) \cap \Gamma_3(v) \cap \Gamma_3(w) = \emptyset$ because none of the $u,v,w$ are on the same grid lines;
b) if exactly one of $U,V,W$ equals 3 then $U = 3$ and $V,W \neq 3$ so $v,w$ are on the same grid line $L$ but $u,v$ and $u,w$ are not. There is exactly one grid line through $u$ that intersects L and that is at a vertex $z$;
c) if exactly two of $U,V,W$ equals 3 then $U,V = 3$ and $W \neq 3$ so $u,w$ and $v,w$ are on grid lines $L$ and $M$ that meet at $w$ but $u,v$ is not. Thus any vertex $z$ that is on a grid line with $u$ and $v$ is diagonally opposite of $w$ so cannot share a grid line with it;
d) if $U = V = W = 3$ then $u,v,w$ are all on the same grid line $L$ so there are 53 vertices on $L$ are the only ones that share this grid line.

*Lemma 3.* Sporadic constraints:

a) $x(1,3,3,3,2,2) = 1$;
b) $x(1,3,3,B), x(2,3,3,B), x(3,1,3,B), x(3,2,3,B), x(3,3,1,B), x(3,3,2,B) \leq 2$ where $B = 2,2,2$.

*Proof.* Again we work with the grid lines. Then

a) if $U = 3$ and $V, W = 2$ then $v, w$ are on the same grid line $L$ but $u, v$ and $u, w$ are not. There is exactly one grid line through $u$ intersecting $L$ at some vertex $t$ which is distance 3 from all 3 vertices $u, v, w$. There is a unique grid line $M$ containing $t$ and $u$. Now $u$ is on a second grid line $N$ which has no intersection with $L$. There is a matching between $L$ and $N$ so $u$ is adjacent to exactly one $z$ on $L$. Thus $x(1,3,3,3,2,2) = 1$.
b) if $U = V = W = 2$ then no two of $u, v, w$ are on the same grid line. Given any two of them, say $u$ and $v$, then let the grid lines determined by $u$ be $P$ and $Q$ and the grid lines determined by $v$ be $L$ and $M$. The grid lines $P$ and $M$ intersect at a unique vertex $z_1$ and the grid lines $L$ and $Q$ meet at a different vertex $z_2$. These two vertices are the unique pair of vertices on grid lines with both $u$ and $v$, that is, $d(u, z_i) = d(v, z_i) = 3$. The value $d(w, z_i)$ is completely determined and may rule out one or both of the $z_i$ as solutions. This proves the inequalities $x(3,3,1,B), x(3,3,2,B) \leq 2$. The other four have analogous proofs.

**Particular solutions to all the blocks.** Now we need a particular solution to each of the chosen blocks. From a computer program, we get this set.

B211. (0, 0, 0, 0, 53, 0, 0, 0, 0, 0, 52, 2, 52, 2652, 104, 2, 104, 2, 0, 0, 0, 0, 106, 2, 0, 2, 0)
B221. (0, 0, 0, 1, 51, 2, 0, 0, 0, 1, 51, 2, 51, 2654, 102, 0, 106, 2, 0, 0, 0, 0, 106, 2, 2, 0, 0)
B222. (1, 0, 0, 0, 52, 0, 0, 0, 2, 0, 52, 0, 52, 2652, 106, 0, 106, 0, 0, 0, 2, 0, 106, 0, 2, 0, 0)
B321. (0, 0, 0, 1, 52, 1, 0, 0, 0, 0, 53, 1, 52, 2704, 52, 0, 54, 53, 0, 0, 0, 1, 106, 1, 1, 0, 1)
B322. (0, 1, 0, 1, 50, 1, 0, 1, 1, 0, 52, 0, 52, 2706, 53, 0, 53, 52, 0, 1, 1, 1, 105, 0, 1, 0, 1)
B332. (0, 1, 0, 0, 52, 0, 0, 1, 1, 0, 52, 0, 54, 2757, 0, 0, 53, 53, 0, 1, 1, 0, 53, 53, 1, 0, 0)
B331. (0, 0, 0, 0, 54, 0, 0, 0, 0, 0, 54, 0, 54, 2754, 0, 0, 54, 54, 0, 0, 0, 0, 54, 54, 1, 0, 0)
B333. (0, 0, 0, 0, 53, 1, 0, 1, 0, 0, 53, 1, 53, 2756, 53, 1, 53, 0, 0, 1, 0, 1, 53, 0, 0, 0, 53)

**Constrained solutions.** Here we develop all solutions to the chosen blocks that satisfy all the known constraints. Essentially, we carry out the integer program using the special form of the equations. (These solutions were initially obtained by the computer. We give a complete analysis here as an independent check.) We go through the chosen blocks and show how to find all solutions given the known constraints. To accomplish this, we know solutions look like a particular solution plus a linear combination of the null vectors $n$.

Let $C$ be the matrix of null basis vectors and $x_p$ be a particular solution. Then a general solution has the form

$$x = Cn + x_p \geq 0$$

plus whatever other constraints we have uncovered. The matrix $C$ has some interesting recursive symmetries. It has the form given by

$$A = \begin{bmatrix} 1 & 0 \\ 0 & 1 \\ -1 & -1 \end{bmatrix}, \quad B = \begin{bmatrix} A & 0 \\ 0 & A \\ -A & -A \end{bmatrix}, \quad C = \begin{bmatrix} B & 0 \\ 0 & B \\ -B & -B \end{bmatrix}.$$

So if we let $n = \begin{bmatrix} a & b & c & d & a' & b' & c' & d' \end{bmatrix}^T$, then constraints come in groups of 9:

1) $a \geq -x_p(1)$      10) $a' \geq -x_p(10)$
2) $b \geq -x_p(2)$      11) $b' \geq -x_p(11)$
3) $a + b \leq x_p(3)$      12) $a' + b' \leq x_p(12)$
4) $c \geq -x_p(4)$      13) $c' \geq -x_p(13)$
5) $d \geq -x_p(5)$      14) $d' \geq -x_p(14)$
6) $c + d \leq x_p(6)$      15) $c' + d' \leq x_p(15)$
7) $a + c \leq x_p(7)$      16) $a' + c' \leq x_p(16)$
8) $b + d \leq x_p(8)$      17) $b' + d' \leq x_p(17)$
9) $a + b + c + d \geq -x_p(9)$      18) $a' + b' + c' + d' \geq -x_p(18)$

19) $a + a' \leq x_p(19)$
20) $b + b' \leq x_p(20)$
21) $a + b + a' + b' \geq -x_p(21)$
22) $c + c' \leq x_p(22)$
23) $d + d' \leq x_p(23)$
24) $c + d + c' + d' \geq -x_p(24)$
25) $a + c + a' + c' \geq -x_p(25)$
26) $b + d + b' + d' \geq -x_p(26)$
27) $a + b + c + d + a' + b' + c' + d' \leq x_p(27)$

.

In addition, from Lemma 2, for every block $x(27) = x_p(27)$ This means that we always have the constraint

Z) $a + b + c + d + a' + b' + c' + d' = 0$

which is stronger than E27.

**Block 333**. To illustrate how this works, we use the particular solution for block 333 and conclude there is but one solution:

$$x_p = [\,0 \quad 0 \quad 0 \quad 0 \quad 53 \quad 1 \quad 0 \quad 1 \quad 0 \quad 0 \quad 53 \quad 1 \quad 53 \quad 2756$$

$$53 \quad 1 \quad 53 \quad 0 \quad 0 \quad 1 \quad 0 \quad 1 \quad 53 \quad 0 \quad 0 \quad 0 \quad 53\,]$$

Clearly coordinate 14 is not helpful. But constraint Z

$$a+b+c+d+a'+b'+c'+d'=0$$

is very strong.

From the 13 zero entries 1, 2, 3, 4, 7, 9, 10, 18, 19, 21, 24, 25, 26 we are now able to prove that $n=0$. From E1, E2 and E3, we get

$$a=b=0.$$

From E4 and E7 we get

$$c=0.$$

From E9 we get

$$d \geq 0.$$

From E18 we get

$$a'+b'+c'+d' \geq 0.$$

Now put this into the equality Z to get

$$0 = a+b+c+d+a'+b'+c'+d' = d+(a'+b'+c'+d')$$

which yields

$$d = 0$$

$$a'+b'+c'+d' = 0.$$

But E21 and E24 give

$$a'+b' \geq 0$$

$$c'+d' \geq 0$$

which yields

$$a'+b' = 0$$

$$c'+d' = 0$$

while E25 and E26 give

$a' + c' \geq 0$

$b' + d' \geq 0$

so

$a' + c' = 0$

$b' + d' = 0$.

These last 4 equations show that the prime variables are also all zero. The only null vector solution is

$$\boxed{0\ 0\ 0\ 0\ 0\ 0\ 0\ 0}.$$

**Block 211**. In this block, we are given $x(27) = x_p(27) = 0$ where

$$x_p = [\,0\quad 0\quad 0\quad 0\quad 53\quad 0\quad 0\quad 0\quad 0\quad 0\quad 52\quad 2\quad 52\quad 2652$$

$$104\quad 2\quad 104\quad 2\quad 0\quad 0\quad 0\quad 0\quad 106\quad 2\quad 0\quad 2\quad 0\,].$$

Again, we have a lot of zeroes: indices 1, 2, 3, 4, 6, 7, 8, 9, 10, 19, 20, 21, 22, 25, 27.

We start with the constraint Z

$a + b + c + d + a' + b' + c' + d' = 0$

and show all the variables are zero.

From E1, E2 and E3, we get

$a = b = 0$.

From E4 and E7 we get

$c = 0$.

From E6 and E9 we get

$d = 0$.

E19, E20 and E21 give us

$$a' = b' = 0.$$

E22 and E25 produce

$$0 \geq c' \geq 0$$

so that leaves only $d'$ to deal with which is the only variable left in the Z. This shows $n = 0$ and there is only one solution. Again, the solution is unique

$$\boxed{0\ 0\ 0\ 0\ 0\ 0\ 0\ 0}.$$

**Block 221**. In this block, we are given $x(27) = x_p(27) = 0$ where

$$x_p = [\,0\ \ 0\ \ 0\ \ 1\ \ 51\ \ 2\ \ 0\ \ 0\ \ 0\ \ 1\ \ 51\ \ 2\ \ 51\ \ 2654$$

$$102\ \ 0\ \ 106\ \ 2\ \ 0\ \ 0\ \ 0\ \ 0\ \ 106\ \ 2\ \ 2\ \ 0\ \ 0\,].$$

Again, we have a lot of zeroes: indices 1, 2, 3, 7, 8, 9, 16, 19, 20, 21, 22, 26, 27.

We start with the constraint Z

$$a + b + c + d + a' + b' + c' + d' = 0.$$

From E1, E2 and E3, we get

$$a = b = 0.$$

E7, E8 and E9 give us

$$c \leq 0$$

$$d \leq 0$$

$$c + d \geq 0$$

and so

$$c = d = 0.$$

E16 then gives

$a' + c' \leq 0$.

E19, E20 and E21 give

$a' = b' = 0$.

E22 gives

$c' \leq 0$.

E26 give

$d' \geq 0$

and in fact from equality Z

$c' = -d'$.

At this point, we have exhausted the zero entries and we have to move on to some other small ones that constrain the two remaining coefficients. We have E25 which yields

$c' \geq -2$

and in fact, each of the possible pairs of values

$c' = d' = 0$

$c' = -d' = -1$

$c' = -d' = -2$

produces a possible solution. The three solutions are

$$\begin{vmatrix} 0 & 0 & 0 & 0 & 0 & 0 & 0 & 0 \\ 0 & 0 & 0 & 0 & 0 & 0 & -1 & 1 \\ 0 & 0 & 0 & 0 & 0 & 0 & -2 & 2 \end{vmatrix}.$$

**Block 321**. In this block, we are given $x(27) = x_p(27) = 1$ where

$$x_p = [\, 0 \quad 0 \quad 0 \quad 1 \quad 52 \quad 1 \quad 0 \quad 0 \quad 0 \quad 0 \quad 53 \quad 1 \quad 52 \quad 2704$$

52  0  54  53  0  0  0  1  106  1  1  0  1 ].

Again, we have a lot of zeroes: indices 1, 2, 3, 7, 8, 9, 10, 16, 19, 20, 21, 26.

We start with the constraint Z

$$a+b+c+d+a'+b'+c'+d'=0.$$

From E1, E2 and E3, we get

$$a=b=0.$$

E7, E8 and E9 give us

$$c \leq 0$$

$$d \leq 0$$

$$c+d \geq 0$$

and so

$$c=d=0.$$

E16 gives

$$c' \leq 0.$$

E19, E20 and E21 give

$$a'=b'=0.$$

E26 gives

$$d' \geq 0.$$

Finally, we employ Z to get

$$c'=-d'.$$

We need to resolve the possible values for $c', d'$. E25 gives us what we need

$$c' \geq -1.$$

Each of the following pairs produces a possible solution

$c' = d' = 0$

$c' = -d' = -1$.

The two solutions are:

$$\begin{vmatrix} 0 & 0 & 0 & 0 & 0 & 0 & 0 & 0 \\ 0 & 0 & 0 & 0 & 0 & 0 & -1 & 1 \end{vmatrix}.$$

**Block 331.** In this block, we are given $x(27) = x_p(27) = 0$ where

$$x_p = [\,0 \quad 0 \quad 0 \quad 0 \quad 54 \quad 0 \quad 0 \quad 0 \quad 0 \quad 0 \quad 54 \quad 0 \quad 54 \quad 2754$$

$$0 \quad 0 \quad 54 \quad 54 \quad 0 \quad 0 \quad 0 \quad 0 \quad 54 \quad 54 \quad 1 \quad 0 \quad 0\,].$$

Again, we have zeroes: indices 1, 2, 3, 4, 6, 7, 8, 9, 10, 12, 15, 16, 19, 20, 21, 22, 26. But we are also going to need the entry with value 1: 25.

We utilize the proven constraint

Z) $a+b+c+d+a'+b'+c'+d' = 0$.

Also

E1) $a \geq 0$
E2) $b \geq 0$
E3) $a+b \leq 0$

so

$a = b = 0$.

Then

E4) $c \geq -x_p(4)$
E7) $c \leq 0$

so

$c = 0$.

Then

E6) $d \leq 0$
E9) $d \geq 0$

so we have

$a = b = c = d = 0$.

Then

E10) $a' \geq 0$
E19) $a' \leq 0$

so

$a' = 0$.

Then

E12) $b' \leq 0$
E21) $b' \geq 0$

so

$b' = 0$.

From Z) we get

$c' + d' = 0$.

Finally

E25) $c' \geq -1$
E22) $c' \leq 0$

so each of the following pairs produces a possible solution

$c' = d' = 0$

$c' = -d' = -1$.

The two solutions are:

$$\begin{vmatrix} 0 & 0 & 0 & 0 & 0 & 0 & 0 & 0 \\ 0 & 0 & 0 & 0 & 0 & 0 & -1 & 1 \end{vmatrix}.$$

**Block 322.** In this block, we are given $x(27) = x_p(27) = 1$ and from Lemma 3, $1 = x(9) = x_p(9) + n(9)$ where

$$x_p = [\,0 \quad 1 \quad 0 \quad 1 \quad 50 \quad 1 \quad 0 \quad 1 \quad 1 \quad 0 \quad 52 \quad 0 \quad 52 \quad 2706$$

$$53 \quad 0 \quad 53 \quad 52 \quad 0 \quad 1 \quad 1 \quad 1 \quad 105 \quad 0 \quad 1 \quad 0 \quad 1\,].$$

Again, we have zeroes: indices 1, 3, 7, 10, 12, 16, 19, 24, 21, 26. But we are also going to need some entries with value 1: 2, 4, 6, 8, 9, 20, 21, 22, 25.

We note the constraints Z

$$a+b+c+d+a'+b'+c'+d' = 0$$

and N

$$a+b+c+d = 0.$$

From E1, E10 and E19 we get

$$a \geq 0$$

$$a' \geq 0$$
$$a+a' \leq 0$$

which implies

$$a = a' = 0,$$

By E3, E12 and E21 we have

$$b \leq 0$$

$$b' \leq 0$$

$$-1 \leq b+b' \leq 0$$

which gives three possible pairs for $(b,b')$

$(0,0), (0,-1), (-1,0)$.

From E7, E16 and E25 we have

$c \leq 0$

$c' \leq 0$

$-1 \leq c + c' \leq 0$

which gives three possible pairs for $(c, c')$

$(0,0), (0,-1), (-1,0)$.

Now from N, we get

$d = -b - c$

and from Z, we get

$d' = -b' - c'$.

This gives us 9 solutions depending on the solutions for $(b, b', c, c')$ in

$\{(0,0), (0,-1), (-1,0)\} \times \{(0,0), (0,-1), (-1,0)\}$.

The nine solutions are

$$\begin{vmatrix} 0 & 0 & 0 & 0 & 0 & 0 & 0 & 0 \\ 0 & 0 & 0 & 0 & 0 & 0 & -1 & 1 \\ 0 & 0 & -1 & 1 & 0 & 0 & 0 & 0 \\ 0 & 0 & 0 & -1 & 0 & -1 & 0 & 0 \\ 0 & 0 & 0 & -1 & 0 & -1 & -1 & -1 \\ 0 & 0 & -1 & -1 & 0 & -1 & 0 & -1 \\ 0 & -1 & 0 & -1 & 0 & 0 & 0 & 0 \\ 0 & -1 & 0 & -1 & 0 & 0 & -1 & -1 \\ 0 & -1 & -1 & -1 & 0 & 0 & 0 & -1 \end{vmatrix}.$$

**Block 222.** In this block, we need $x(9)$, $x(18)$, $x(21)$, $x(24)$, $x(25)$, $x(26) \leq 2$ and the equality Z, $x(27) = x_p(27) = 0$ where

$$x_p = [\ 1\quad 0\quad 0\quad 0\quad 52\quad 0\quad 0\quad 0\quad 2\quad 0\quad 52\quad 0\quad 52\quad 2652$$

$$106\quad 0\quad 106\quad 0\quad 0\quad 0\quad 2\quad 0\quad 106\quad 0\quad 2\quad 0\quad 0\ ].$$

Again, we have zeroes: indices 2, 3, 4, 6, 7, 8, 10, 12, 16, 18, 19, 20, 22, 24, 26, 27 and also $x_p(1) = 1$, $x_p(9) = 2$, $x_p(21) = 2$, $x_p(25) = 2$.

We note the constraint

Z) $a + b + c + d + a' + b' + c' + d' = 0$.

We use the zeroes and the one to rewrite the matrix constraints:

E1) $a \geq -1$

E2) $b \geq 0$

E3) $a + b \leq 0$

E4) $c \geq 0$

E6) $c + d \leq 0$

E7) $a + c \leq 0$

E8) $b + d \leq 0$

E10) $a' \geq 0$

E12) $a' + b' \leq 0$

E16) $a' + c' \leq 0$

E18) $a' + b' + c' + d' \geq 0$

E19) $a + a' \leq 0$

E20) $b + b' \leq 0$

E22) $c + c' \leq 0$

E24) $c + c' + d + d' \geq 0$

E26) $b+b'+d+d' \geq 0$.

The six constraints from Lemma 3, $x(9), x(18), x(21), x(24), x(25), x(26) \leq 2$, will be useful. We record those here:

E9) $-2 \leq a+b+c+d \leq 0$

E18) $0 \leq a'+b'+c'+d' \leq 2$

E21) $-2 \leq a+b+a'+b' \leq 0$

E24) $0 \leq c+d+c'+d' \leq 2$

E25) $-2 \leq a+a'+c+c' \leq 0$

E26) $0 \leq b+b'+d+d' \leq 2$.

First we use E1, E2 and E3 to get

$a+1 \geq 0$
$b \geq 0$
$0 \leq (a+1)+b \leq 1$

which yields three possible choices for $(a,b)$: $(a,b) = (-1,0)$, $(-1,1)$, $(0,0)$.

Then we use E1, E4 and E7 in identical fashion to get $(a,c) = (-1,0), (-1,1), (0,0)$.

Finally, E1, E10 and E19 also have the same relations so $(a,a') = (-1,0)$, $(-1,1), (0,0)$.

Put this all together to get choices for $(a,b,c,a')$:

1) $(a,b,c,a') = (0,0,0,0)$
2) $(a,b,c,a') = (-1,0,0,0)$
3) $(a,b,c,a') = (-1,0,0,1)$
4) $(a,b,c,a') = (-1,0,1,0)$
5) $(a,b,c,a') = (-1,0,1,1)$
6) $(a,b,c,a') = (-1,1,0,0)$
7) $(a,b,c,a') = (-1,1,0,1)$
8) $(a,b,c,a') = (-1,1,1,0)$
9) $(a,b,c,a') = (-1,1,1,1)$.

We can step through them one at a time. But first we make explicit the method we use to resolve each case.

First

E6) and E8) $d \leq \min(-b, -c)$
E12) and E20) $b' \leq \min(-b, -a')$
E16) and E22) $c' \leq \min(-c, -a')$.
E9) $-2 - (a + b + c) \leq d \leq -(a + b + c)$
E21) $-2 - (a + b + a') \leq b' \leq -(a + b + a')$
E25) $-2 - (a + c + a') \leq c' \leq -(a + c + a')$

We can rewrite these as

$-2 - (a + b + c) \leq d \leq \min(-(a + b + c), -b, -c)$
$-2 - (a + b + a') \leq b' \leq \min(-(a + b + a'), -b, -a')$
$-2 - (a + c + a') \leq c' \leq \min(-(a + c + a'), -c, -a')$

We also need

Z) $d' = -(a + b + c + a') - (d + b' + c')$.

Case 1. $(a, b, c, a') = (0, 0, 0, 0)$

$-2 \leq d \leq 0$
$-2 \leq b' \leq 0$
$-2 \leq c' \leq 0$
$d' = -(b' + c' + d)$

All 27 of these possibilities are allowed.

Case 2. $(a, b, c, a') = (-1, 0, 0, 0)$

$-1 \leq d \leq 0$
$-1 \leq b' \leq 0$
$-1 \leq c' \leq 0$
$d' = 1 - (b' + c' + d)$.

All 8 of these possibilities are allowed.

Case 3. $(a, b, c, a') = (-1, 0, 0, 1)$

$-1 \leq d \leq 0$
$-2 \leq b' \leq -1$
$-2 \leq c' \leq -1$
$d' = -(b' + c' + d)$.

All 8 of these possibilities are allowed.

Case 4. $(a, b, c, a') = (-1, 0, 1, 0)$

$-2 \leq d \leq -1$
$-1 \leq b' \leq 0$
$-2 \leq c' \leq -1$
$d' = -(b' + c' + d)$.

All 8 of these possibilities are allowed.

Case 5. $(a, b, c, a') = (-1, 0, 1, 1)$

$-2 \leq d \leq -1$
$-2 \leq b' \leq -1$
$-3 \leq c' \leq -1$
$d' = -1 - (b' + c' + d)$.

All 12 of these possibilities are allowed.

Case 6. $(a, b, c, a') = (-1, 1, 0, 0)$

$-2 \leq d \leq -1$
$-2 \leq b' \leq -1$
$-1 \leq c' \leq 0$
$d' = -(b' + c' + d)$.

All 8 of these possibilities are allowed.

Case 7. $(a, b, c, a') = (-1, 1, 0, 1)$

$-2 \leq d \leq -1$
$-3 \leq b' \leq -1$
$-2 \leq c' \leq -1$
$d' = -1 - (b' + c' + d)$.

All 12 of these possibilities are allowed.

Case 8. $(a, b, c, a') = (-1, 1, 1, 0)$

$-3 \leq d \leq -1$
$-2 \leq b' \leq -1$
$-2 \leq c' \leq -1$
$d' = -1 - (b' + c' + d)$.

All 12 of these possibilities are allowed.

Case 9. $(a, b, c, a') = (-1, 1, 1, 1)$

$-3 \leq d \leq -1$
$-3 \leq b' \leq -1$
$-3 \leq c' \leq -1$
$d' = -2 - (b' + c' + d)$.

All 27 of these possibilities are allowed. We can enumerate the 122 solutions in case order:

Case 1

| | | | | | | | |
|---|---|---|---|---|---|---|---|
| 0 | 0 | 0 | 0 | −2 | −2 | −2 | 6 |
| 0 | 0 | 0 | 0 | −2 | −2 | −1 | 5 |
| 0 | 0 | 0 | 0 | −2 | −2 | 0 | 4 |
| 0 | 0 | 0 | 0 | −2 | −1 | −2 | 5 |
| 0 | 0 | 0 | 0 | −2 | −1 | −1 | 4 |
| 0 | 0 | 0 | 0 | −2 | −1 | 0 | 3 |
| 0 | 0 | 0 | 0 | −2 | 0 | −2 | 4 |
| 0 | 0 | 0 | 0 | −2 | 0 | −1 | 3 |
| 0 | 0 | 0 | 0 | −2 | 0 | 0 | 2 |
| 0 | 0 | 0 | 0 | −1 | −2 | −2 | 5 |
| 0 | 0 | 0 | 0 | −1 | −2 | −1 | 4 |
| 0 | 0 | 0 | 0 | −1 | −2 | 0 | 3 |
| 0 | 0 | 0 | 0 | −1 | −1 | −2 | 4 |
| 0 | 0 | 0 | 0 | −1 | −1 | −1 | 3 |
| 0 | 0 | 0 | 0 | −1 | −1 | 0 | 2 |
| 0 | 0 | 0 | 0 | −1 | 0 | −2 | 3 |
| 0 | 0 | 0 | 0 | −1 | 0 | −1 | 2 |
| 0 | 0 | 0 | 0 | −1 | 0 | 0 | 1 |
| 0 | 0 | 0 | 0 | 0 | −2 | −2 | 4 |
| 0 | 0 | 0 | 0 | 0 | −2 | −1 | 3 |
| 0 | 0 | 0 | 0 | 0 | −2 | 0 | 2 |
| 0 | 0 | 0 | 0 | 0 | −1 | −2 | 3 |
| 0 | 0 | 0 | 0 | 0 | −1 | −1 | 2 |
| 0 | 0 | 0 | 0 | 0 | −1 | 0 | 1 |
| 0 | 0 | 0 | 0 | 0 | 0 | −2 | 2 |
| 0 | 0 | 0 | 0 | 0 | 0 | −1 | 1 |
| 0 | 0 | 0 | 0 | 0 | 0 | 0 | 0 |

Case 2

| -1 | 0 | 0 | 0 | -1 | -1 | -1 | 4 |
| -1 | 0 | 0 | 0 | -1 | -1 | 0 | 3 |
| -1 | 0 | 0 | 0 | -1 | 0 | -1 | 3 |
| -1 | 0 | 0 | 0 | -1 | 0 | 0 | 2 |
| -1 | 0 | 0 | 0 | 0 | -1 | -1 | 3 |
| -1 | 0 | 0 | 0 | 0 | -1 | 0 | 2 |
| -1 | 0 | 0 | 0 | 0 | 0 | -1 | 2 |
| -1 | 0 | 0 | 0 | 0 | 0 | 0 | 1 |

Case 3

| -1 | 0 | 0 | 1 | -1 | -2 | -2 | 5 |
| -1 | 0 | 0 | 1 | -1 | -2 | -1 | 4 |
| -1 | 0 | 0 | 1 | -1 | -1 | -2 | 4 |
| -1 | 0 | 0 | 1 | -1 | -1 | -1 | 3 |
| -1 | 0 | 0 | 1 | 0 | -2 | -2 | 4 |
| -1 | 0 | 0 | 1 | 0 | -2 | -1 | 3 |
| -1 | 0 | 0 | 1 | 0 | -1 | -2 | 3 |
| -1 | 0 | 0 | 1 | 0 | -1 | -1 | 2 |

Case 4

| -1 | 0 | 1 | 0 | -2 | -1 | -2 | 5 |
| -1 | 0 | 1 | 0 | -2 | -1 | -1 | 4 |
| -1 | 0 | 1 | 0 | -2 | 0 | -2 | 4 |
| -1 | 0 | 1 | 0 | -2 | 0 | -1 | 3 |
| -1 | 0 | 1 | 0 | -1 | -1 | -2 | 4 |
| -1 | 0 | 1 | 0 | -1 | -1 | -1 | 3 |
| -1 | 0 | 1 | 0 | -1 | 0 | -2 | 3 |
| -1 | 0 | 1 | 0 | -1 | 0 | -1 | 2 |

Case 5

| -1 | 0 | 1 | 1 | -2 | -2 | -3 | 6 |
| -1 | 0 | 1 | 1 | -2 | -2 | -2 | 5 |
| -1 | 0 | 1 | 1 | -2 | -2 | -1 | 4 |
| -1 | 0 | 1 | 1 | -2 | -1 | -3 | 5 |
| -1 | 0 | 1 | 1 | -2 | -1 | -2 | 4 |
| -1 | 0 | 1 | 1 | -2 | -1 | -1 | 3 |
| -1 | 0 | 1 | 1 | -1 | -2 | -3 | 5 |
| -1 | 0 | 1 | 1 | -1 | -2 | -2 | 4 |
| -1 | 0 | 1 | 1 | -1 | -2 | -1 | 3 |
| -1 | 0 | 1 | 1 | -1 | -1 | -3 | 4 |
| -1 | 0 | 1 | 1 | -1 | -1 | -2 | 3 |
| -1 | 0 | 1 | 1 | -1 | -1 | -1 | 2 |

Case 6

| -1 | 1 | 0 | 0 | -2 | -2 | -1 | 5 |
| -1 | 1 | 0 | 0 | -2 | -2 | 0 | 4 |
| -1 | 1 | 0 | 0 | -2 | -1 | -1 | 4 |
| -1 | 1 | 0 | 0 | -2 | -1 | 0 | 3 |
| -1 | 1 | 0 | 0 | -1 | -2 | -1 | 4 |
| -1 | 1 | 0 | 0 | -1 | -2 | 0 | 3 |
| -1 | 1 | 0 | 0 | -1 | -1 | -1 | 3 |
| -1 | 1 | 0 | 0 | -1 | -1 | 0 | 2 |

Case 7

| -1 | 1 | 0 | 1 | -2 | -3 | -2 | 6 |
| -1 | 1 | 0 | 1 | -2 | -3 | -1 | 5 |
| -1 | 1 | 0 | 1 | -2 | -2 | -2 | 5 |
| -1 | 1 | 0 | 1 | -2 | -2 | -1 | 4 |
| -1 | 1 | 0 | 1 | -2 | -1 | -2 | 4 |
| -1 | 1 | 0 | 1 | -2 | -1 | -1 | 3 |
| -1 | 1 | 0 | 1 | -1 | -3 | -2 | 5 |
| -1 | 1 | 0 | 1 | -1 | -3 | -1 | 4 |
| -1 | 1 | 0 | 1 | -1 | -2 | -2 | 4 |
| -1 | 1 | 0 | 1 | -1 | -2 | -1 | 3 |
| -1 | 1 | 0 | 1 | -1 | -1 | -2 | 3 |
| -1 | 1 | 0 | 1 | -1 | -1 | -1 | 2 |

Case 8

| -1 | 1 | 1 | 0 | -3 | -2 | -2 | 6 |
| -1 | 1 | 1 | 0 | -3 | -2 | -1 | 5 |
| -1 | 1 | 1 | 0 | -3 | -1 | -2 | 5 |
| -1 | 1 | 1 | 0 | -3 | -1 | -1 | 4 |
| -1 | 1 | 1 | 0 | -2 | -2 | -2 | 5 |
| -1 | 1 | 1 | 0 | -2 | -2 | -1 | 4 |
| -1 | 1 | 1 | 0 | -2 | -1 | -2 | 4 |
| -1 | 1 | 1 | 0 | -2 | -1 | -1 | 3 |
| -1 | 1 | 1 | 0 | -1 | -2 | -2 | 4 |
| -1 | 1 | 1 | 0 | -1 | -2 | -1 | 3 |
| -1 | 1 | 1 | 0 | -1 | -1 | -2 | 3 |
| -1 | 1 | 1 | 0 | -1 | -1 | -1 | 2 |

Case 9

$$\begin{bmatrix}
-1 & 1 & 1 & 1 & -3 & -3 & -3 & 7 \\
-1 & 1 & 1 & 1 & -3 & -3 & -2 & 6 \\
-1 & 1 & 1 & 1 & -3 & -3 & -1 & 5 \\
-1 & 1 & 1 & 1 & -3 & -2 & -3 & 6 \\
-1 & 1 & 1 & 1 & -3 & -2 & -2 & 5 \\
-1 & 1 & 1 & 1 & -3 & -2 & -1 & 4 \\
-1 & 1 & 1 & 1 & -3 & -1 & -3 & 5 \\
-1 & 1 & 1 & 1 & -3 & -1 & -2 & 4 \\
-1 & 1 & 1 & 1 & -3 & -1 & -1 & 3 \\
-1 & 1 & 1 & 1 & -2 & -3 & -3 & 6 \\
-1 & 1 & 1 & 1 & -2 & -3 & -2 & 5 \\
-1 & 1 & 1 & 1 & -2 & -3 & -1 & 4 \\
-1 & 1 & 1 & 1 & -2 & -2 & -3 & 5 \\
-1 & 1 & 1 & 1 & -2 & -2 & -2 & 4 \\
-1 & 1 & 1 & 1 & -2 & -2 & -1 & 3 \\
-1 & 1 & 1 & 1 & -2 & -1 & -3 & 4 \\
-1 & 1 & 1 & 1 & -2 & -1 & -2 & 3 \\
-1 & 1 & 1 & 1 & -2 & -1 & -1 & 2 \\
-1 & 1 & 1 & 1 & -1 & -3 & -3 & 5 \\
-1 & 1 & 1 & 1 & -1 & -3 & -2 & 4 \\
-1 & 1 & 1 & 1 & -1 & -3 & -1 & 3 \\
-1 & 1 & 1 & 1 & -1 & -2 & -3 & 4 \\
-1 & 1 & 1 & 1 & -1 & -2 & -2 & 3 \\
-1 & 1 & 1 & 1 & -1 & -2 & -1 & 2 \\
-1 & 1 & 1 & 1 & -1 & -1 & -3 & 3 \\
-1 & 1 & 1 & 1 & -1 & -1 & -2 & 2 \\
-1 & 1 & 1 & 1 & -1 & -1 & -1 & 1
\end{bmatrix}.$$

**Block 332.** In this block, we are given $x(27) = x_p(27) = 0$ where

$$x_p = [\,0 \quad 1 \quad 0 \quad 0 \quad 52 \quad 0 \quad 0 \quad 1 \quad 1 \quad 0 \quad 52 \quad 0 \quad 54 \quad 2757$$
$$0 \quad 0 \quad 53 \quad 53 \quad 0 \quad 1 \quad 1 \quad 0 \quad 53 \quad 53 \quad 1 \quad 0 \quad 0\,].$$

Again, we have zeroes: indices 1, 3, 4, 6, 7, 10, 12, 15, 16, 19, 22, 26. But we are also going to need some entries with value 1: 2, 8, 9, 20, 21, 25.

Again we utilize

Z) $a+b+c+d+a'+b'+c'+d'=0$.

Then

E1) $a \geq 0$
E10) $a' \geq 0$
E19) $a+a' \leq 0$

so

$a = a' = 0$.

E4) $c \geq 0$
E7) $c \leq 0$

so

$c = 0$.

Now

E3) $b \leq 0$
E6) $d \leq 0$

combines with E9 to give

$0 \geq b+d \geq -1$.

We also have

E12) $b' \leq 0$
E15) $c'+d' \leq 0$
Z) $b+d+b'+c'+d' = 0$

so

$b' = 0$
$c' + d' = 0$
$b + d = 0$

but then

$b = d = 0$.

At this point we have

$a = a' = c = b' = b = d = 0$.

Now

E16) $c' \leq 0$
E25) $c' \geq -1$.

Each of the following pairs produces a possible solution

$c' = d' = 0$
$c' = -d' = -1$.

The two solutions are

$$\begin{vmatrix} 0 & 0 & 0 & 0 & 0 & 0 & 0 & 0 \\ 0 & 0 & 0 & 0 & 0 & 0 & -1 & 1 \end{vmatrix}.$$

To summarize the number of solutions for each block:

| Block | 333 | 211 | 221 | 321 | 331 | 322 | 222 | 332 |
|-------|-----|-----|-----|-----|-----|-----|-----|-----|
| Count | 1   | 1   | 3   | 2   | 2   | 9   | 122 | 2   |

**Discussion**. Lemma 4 in [1] uses some side conditions to show that $x(333221) = 0$ and in fact all of the solutions to Block 221 have that constraint. Note that $x(132221) = 0$ because of the triangle inequality. Note that our three different solutions for Block 221 have three different values for $x(331221)$ ranging from 0 to 2. So it is not clear to us as to why Lemma 4 in [1] claims the value is 2. The very last paragraph of the paper [1] claims a contradiction to the existence assumption because the value of $x(331221)$ cannot be 2 and in fact, we have shown it does not have to be. The paper claims that $x(221221) = 51$ but we do not understand why this is so. The value 51 occurs with the value $x(331221) = 2$ but it appears that $x(221221) - x(331221) = 49$. If there is some side condition that contradicts that then the existence would be disproven with Block 221 alone and everything else would be irrelevant.

We have shown that although counting features based on the six possible distances between four points to disprove existence of a Moore graph seems like a good idea, it seems to be insufficient because there are many natural number solutions. In fact, we have shown that there exist relationships between these features that depend on the properties of a Moore graph and are independent of its existence. Unfortunately, these relations seem to say nothing about how to construct a Moore graph if it exists. For this reason, we offer an alternative algebraic description that gives a necessary and sufficient condition for existence.

**Necessary and sufficient properties for existence**. We can choose a single vertex $v$ in the graph and try to build the graph by breadth first search. There are $d$ neighbors $N(v)$ and none of them can be adjacent or there would be triangles. The remaining vertices $V$ of $M$ are adjacent to the vertices in $N(v)$. We propose searching for the graph $H$ that forms the induced subgraph of $V$. We show that $H$ exists if and only if there exists a collection of permutations among a certain family of permutations which are fixed-point free.

*Theorem 4.* The Moore Graph $M$ exists if and only if $H$ exists with the following properties:

i) $H$ is $d$-partite;
ii) each part has $d-1$ vertices;
iii) the degree of each vertex in $H$ is $d-1$;
iv) each vertex in one part is adjacent to exactly one vertex in any other part;
v) $H$ has no triangles or squares.

*Proof.* This theorem is just a direct encoding of the properties of $H$. The $d$ parts are just the sets adjacent to a common member of $N(v)$. Given a vertex $x$ in one part, it can't be adjacent to two vertices $u$ and $v$ in another part or else $x$, $u$ and $v$ would be part of a square.

*Corollary 5.* The Moore Graph $M$ exists if and only if there exists $\binom{d}{2}$ permutations $\theta_{ij}$ with $1 \leq i < j \leq d$ on a set with $d-1$ elements so that all of the following permutations are fixed point free where $(i,j,k,l)$ are distinct (with $\theta_{ji} = \theta_{ij}^{-1}$ when $i < j$) :

a) $\theta_{ki}\theta_{ij}\theta_{jk}$ ;

b) $\theta_{lk}\theta_{ki}\theta_{ij}\theta_{jl}$.

*Proof.* Let $X_i$ be the set of vertices in part $i$ of the graph $H$, indexed by $X = \{1, 2, ..., d-1\}$. The edges from $X_i$ to $X_j$ form a bijection, hence a permutation $\theta_{ij}$ on $X$. A fixed point $x$ for one of the listed expressions then represents a vertex on a cycle in $H$ which we want to forbid. Equation *a)* says $x \in X_k$ is involved in a triangle. Equation *b)* says $x \in X_l$ is involved in a square.

It turns out that this lemma is equivalent to the construction in Theorem 1 as the next corollary shows.

*Corollary 6.* The permutations in Corollary 5 exist if and only if there exists $\binom{d-1}{2}$ permutations $\psi_{ij}$ with $1 \leq i < j \leq d-1$ on a set with $d-1$ elements so that all of the following permutations are fixed point free where $(i, j, k, l)$ with $i < j < k < l$ (with $\psi_{ji} = \psi_{ij}^{-1}$ when $i < j$):

a) $\psi_{ij}$;

b) $\psi_{jk}\psi_{ij}$;

c) $\psi_{ki}\psi_{ij}\psi_{jk}$;

d) $\psi_{lk}\psi_{ki}\psi_{ij}\psi_{jl}$.

*Proof.* If the $\psi_{ij}$ are given, we let $\theta_{ij} = \psi_{ij}$ for $j < d$. In addition, define $\theta_{id}$ be the identity permuation. The conditions in Corollary 5 are clearly satisfied as long as no index is $d$. The first two conditions handle those cases. For example, if $\psi_{jk}\psi_{ij}(m) = m$ then we get the contradiction

$$\theta_{id}\theta_{jk}\theta_{ij}\theta_{di}(m) = \psi_{id}\psi_{jk}\psi_{ij}\psi_{di}(m) = \psi_{jk}\psi_{ij}(m) = m.$$

**Questions**. Here are some questions that we haven't been able to answer.

Q1. The given equations provide no relationship between solutions of two different blocks. Are there some more constraints that haven't been considered?

Q2. Assuming that a Moore graph $G$ exists, if $w$ is a vertex how many sets of vertices $x, y, z$ exist so that the collection of 6 distances between the vertices is specified? It seems that the solutions we find fix the distances between $x, y, z$ and count how many $w$ exist with fixed distance to the other three vertices. In order to construct $G$, it seems like it would be useful to fix all the distances and vary the other 3 vertices. This generalizes the notion of the degree of a vertex.

**Final notes**. The intersection numbers are useful to disprove the existence of $\Gamma_3$ and consequently the Moore graph. But since they have solutions, more constraints are needed. However, if we want to construct the Moore graph using the lattice graph, the equations (+) are somewhat irrelevant. If we have a set of permutations as in Corollary 6, then checking that the graph is distance regular is similar to checking if it has girth 5. If it does then this is a Moore graph and properties like (+) are consequences and do not need to be checked.